\documentclass[a4paper,11pt]{article}

\usepackage[pdftex]{graphicx}

\usepackage{graphicx}
\usepackage[pdftex]{color}
\usepackage{amsmath}
\usepackage{amsfonts} 
\usepackage{subfig}
\usepackage{soul}
\usepackage{array}
\usepackage{ragged2e}
\justifying

\newtheorem{propo}{Proposition}[section]
\newtheorem{lemm}{Lemma}[section]

\newtheorem{rem}{Remark}

 \def\diag{\mathrm{diag}}

\def\tr{\mathrm{Tr}}

\newcommand{\tp}{^{\top}}

\newcommand{\beq}{\begin{equation}}
\newcommand{\eeq}{\end{equation}}
\newcommand{\bea}{\begin{eqnarray}}
\newcommand{\eea}{\end{eqnarray}}

\newcommand{\bsea}{\begin{subeqnarray}}
\newcommand{\esea}{\end{subeqnarray}}
\newcommand{\nn}{\nonumber}

\def\bmat{\left[ \begin{array}}
\def\emat{\end{array} \right]}
\def\tr{\text{tr}} 
\def\diag{\text{diag}} 
\def\ofd{\text{ofd}}
\def\argmin{\text{arg\,min}}
\definecolor{Royalblue}{cmyk}{1,0.30,0.2,0.2}

\newcommand{\vale}{\color{black}}
\newcommand{\alg}[1]{\begin{align}#1\end{align}}
 
\newcommand{\proof}{{\em Proof. }} 
\newcommand{\qed}{ $\square$} 

\begin{document}

\title{An alternating minimization algorithm for Factor Analysis}

 \author{Valentina Ciccone, Augusto Ferrante and Mattia Zorzi \footnote{Department of Information Engineering, University of Padova}}

\maketitle

\begin{abstract}
The problem of decomposing a given covariance matrix as the sum of a positive semi-definite matrix of given rank and a positive semi-definite diagonal matrix, is considered.
We present a projection-type algorithm to address this problem. This algorithm appears to perform extremely well and is extremely fast even when the given covariance matrix has a very large dimension. The effectiveness of the algorithm is assessed through simulation studies and by applications to three real datasets that are considered as benchmark for the problem. A local convergence analysis of the algorithm is also presented.
\end{abstract}

\section{Introduction}\label{Intro}

The problem of decomposing a given covariance matrix into the sum of a low rank matrix $L$ plus a diagonal matrix $D$ bursts more than a century of tradition in scientific literature. In fact, it may be viewed as a linear algebraic counterpart of a \textit{Factor Analysis} problem which is a problem in multivariate statistics {\vale aiming to extract statistical commonalities among data.\\ Factor models were first introduced by Spearman \cite{spearman_1904} at the beginning of the last century in the framework of psychological tests. Since then a rich stream of literature has followed combining psychology and mathematics, see for example \cite{BURT_1909}, \cite{Spearman-Holzinger-24}, \cite{kaiser1958varimax}, \cite{tucker1973reliability}, \cite{shapiro1982rank} and \cite{Bekker-deLeeuw_1987} with reference therein. Soon the interest for this type of models has grown significantly also outside the psychological community, see \cite{deistler2007}, \cite{scherrer1998structure},\cite{GEWEKE_DYNAMIC}, \cite{DEISTLER_1997}, \cite{picci1986dynamic}, \cite{bottegal2015modeling} and \cite{zorzi2016ar}, boasting nowadays applications in countless fields of science.
We refer to \cite{finesso2016factor}, \cite{ning2015linear} and \cite{bertsimas2017certifiably} }
for an {\vale up to date} discussion on the importance of the problem, on its applications, on the formidable stream of literature produced on this topic in the last century, and on the numerous variants in which the problem can be formulated.

In this work we take an optimization-oriented viewpoint: for a given covariance matrix $\Sigma$ and a given rank $r$ we want to find a positive semidefinite matrix $L$ with rank not larger than  $r$ and a positive semidefinite diagonal matrix $D$ such that their sum is as closest as possible to $\Sigma$.
A closed-form solution for this problem appears to be out of reach so that numerical techniques are needed. We propose an easy-to-implement iterative algorithm, based on alternating minimization, to solve numerically the considered problem.
This algorithm appears to perform extremely well and in simulations converges very rapidly to the solution.
However, despite the simplicity of the algorithm, the  convergence analysis is non trivial due to the non-convexity of the set of low rank matrices.

The paper is organized as follows. In Section \ref{Section_FA} we first introduce the Factor Analysis problem together with a brief review of the available literature. Then  the  addressed problem is stated and motivated. 
In Section \ref{Section_algorithm} we present the proposed algorithm and  the results of numerical simulations are summarized in Section \ref{sim} together with an application to real datasets.  In Section \ref{Section_convergence} the local convergence analysis for the proposed algorithm is discussed.
Finally, Section \ref{pojection-approach} proposes  a different interpretation of the proposed algorithm that can be viewed as an alternating projection procedure.

\subsection{Notation}
Given a matrix $M$, we denote its transpose by $M\tp$; if $M$ is a square matrix 
$\tr(M)$ denotes its trace. 
The symbol $\mathbf{Q}_n$ denotes the vector space of real symmetric matrices of size $n$. If $X\in\mathbf{Q}_n$ is positive definite or positive semi-definite we write $X\succ 0$ or $X\succeq 0$, respectively.
Moreover, we denote by $\mathbf{D}_n$ the vector space of  diagonal matrices  of size $n$.
We denote by $\ofd(\cdot)$ the self-adjoint operator orthogonally projecting $\mathbf{Q}_n$ onto  the orthogonal complement of $\mathbf{D}_n$ in $\mathbf{Q}_n$, i.e. if $M\in  \mathbf{Q}_n$, $\ofd(M)$ is the matrix in which each off-diagonal element is equal to the corresponding element of $M$ (and each diagonal element is clearly zero).
The Frobenius norm is denoted by $\Vert \cdot \Vert_F$ while $\Vert \cdot \Vert$ denotes the Euclidean norm.

\section{Preliminaries in Factor Analysis and Problem Definition}\label{Section_FA}

Factor models are used to described high dimensional vectors of data in terms of a small number of common latent factors.
In its simplest formulation, the classic (linear static) factor model is given by
 \alg{\label{fact_model} y&=Ax+z}
 where $A\in \mathbb{R}^{n\times r}$, with $r<<n$, is the so-called factor loading matrix, $x$ is the vector of (independent) latent factors and $z$ represents the idiosyncratic component. Here, $x$ and $z$  are zero-mean, independent Gaussian random vectors; the covariance matrix of $x$  is the identity matrix of dimension $r$ and the covariance matrix of $z$ is a diagonal matrix $D\in \mathbf{D}_n$. Note that, $Ax$ represents the latent variable. 
Clearly, $y$ is itself a Gaussian random vector with zero mean and we denote by $\Sigma$ its covariance matrix. Since $x$ and $z$ are independent it holds that
\begin{equation}
\label{decomposition}
\Sigma=L+D
\end{equation}
where $L:=AA^\top$ and $D$ are the covariance matrices of $Ax$ and $z$, respectively. Thus, $L$  has rank equal to $r$, and $D$ is diagonal.

Hence, in its original conception the construction of a factor model is mathematically equivalent to a matrix additive decomposition problem which seeks, for a given $\Sigma$, a decomposition of the type of \eqref{decomposition}.
Of course the model is maximally parsimonious if the rank of $L$ is minimum. 
The problem of minimizing the rank of $L$ in  decomposition (\ref{decomposition})  is known as 
Frisch's problem and, to date, no exact solution for such a problem is actually available, with the only exception of the special case when this minimum rank is $r=n-1$, in which case a closed-form parametrization of the solutions is provided in \cite{reiersol1950identifiability}. This lack of explicit formulas has motivated a rich stream of literature and different numerical approaches which have been proposed over the years.
A relaxation of this problem has also been considered in which the matrix $D$ is only required to be diagonal but not positive semi-definite. This is known as Shapiro's problem.

The main difficulty in these problems is related to the non-convexity of the rank function so that a viable alternative is to consider the so called  \textit{minimum trace factor analysis} problem, \cite{shapiro1982rank}, \cite{MIN_RANK_SHAPIRO_1982}:
\begin{equation}
\label{classic}
\begin{aligned}
\min_{L,D\in\mathbf{Q}_n} \;\; & \tr (L)\\
& L,D\succeq 0 \\
& \Sigma = L+D\\
& D\in \mathbf{D}_n
\end{aligned}
\end{equation}
 where the trace of L is used as convex surrogate of the rank function as shown in \cite{FAZEL_MIN_RANK_APPLICATIONS_2002}, \cite{FAZEL_MINIMUM_RANK_2004}. 
 
Note that, in many cases the equality constraint in \eqref{classic} may be too compelling. Therefore, an alternative approach is to allow for residuals in the decomposition. Typically, this leads to an optimization problem where the residual $\Sigma - L - D$ is minimized with respect to a chosen norm under a constraint limiting the rank of $L$. This approach is known as \textit{minimum residual factor analysis}, see \cite{harman1966factor}, \cite{shapiro2002statistical},\cite{bertsimas2017certifiably}. 
Note that the presence of the rank constraint makes such problems non convex and several heuristic have been proposed to deal with it.\\ 
Other approaches to factor analysis encompass: principal component factor analysis as in \cite{bai2008large}, maximum likelihood methods as in \cite{bai2012statistical}, or the establishing of a certificate of optimal low rank as in \cite{guttman1954some}. 
Moreover, several variants of the mentioned approaches have been proposed by weakening modelling assumptions or by introducing additional constraints for example to account for errors in the covariance matrix estimation as in  \cite{ciccone2017factor}, \cite{ning2015linear} and \cite{agarwal2012noisy}.

The problem we are going to consider is a minimum residual type problem:
for a given $r$ and a given matrix $\Sigma$ we want to find a positive semidefinite matrix $L$ with rank at most $r$ and a positive semidefinite diagonal matrix $D$ such that their sum is as close as possible to $\Sigma$.
This can be formalized as follows:
\begin{equation}
\begin{aligned}
\label{the_problem0}
(L^*,D^*):=\argmin_{L\in \mathcal{L}_{n,r},D\in \mathcal{D}_n} \quad & \Vert \Sigma -L-D \Vert_F^2 
\end{aligned}
\end{equation}
where the sets $\mathcal{D}_n$ and $\mathcal{L}_{n,r}$ are defined as:
\begin{equation*}
\mathcal{L}_{n,r}:=\lbrace X\in\mathbf{Q}_n: X\succeq 0, \; \text{rank}(X) \leq r \rbrace,
\end{equation*}
\begin{equation*}
\mathcal{D}_n:=\lbrace X \in\mathbf{D}_n: X\succeq 0\rbrace.
\end{equation*}
 Note that, in practice, $r$ can be obtained by resorting to available methods for estimating the number of factors as \cite{bai2002determining}, \cite{lam2012factor} and \cite{ciccone2017factor}. Alternatively, the problem can be solved for increasing values of $r$ until the residue $\Vert \Sigma -L^*-D^* \Vert_F$ is not greater than a certain tolerance.  In the case that $\Sigma$ is the sample covariance matrix estimated from the data, this is equivalent to find a good trade-off between the fit term (i.e. the residue) and the complexity of the model (i.e. $r$).

{\vale
Our approach is close in spirit to the one proposed by } \cite{bertsimas2017certifiably} where the q-norm of the residue is minimized under the following ulterior constraint:
\beq\label{additionalconstr}
\Sigma-D \succeq 0.
\eeq
This constraint is perfectly justified if we assume that the covariance matrix  $\Sigma$ of
$y$ is known with great precision, $r$ is the number of the most significant common factors and the residue $\Sigma -L^*-D^* $ accounts for other common factors that are less significant.
Our approach considers instead the case (that is realistic in many practical situations) in which $\Sigma$ has been estimated form the data and is therefore only an approximation of the ``true'' covariance matrix; for more details on this case we refer the reader to  \cite{bai2002determining}, \cite{lam2012factor}, \cite{ciccone2017factor} and \cite{fan2013large}. In this setting, the residue $\Sigma -L^*-D^* $ accounts also for 
the uncertainty in the estimation of $\Sigma$ so that the constraint $\Sigma-D \succeq 0$ must  not 
be imposed.
Of course, if we find an exact decomposition $\Sigma=L^*+D^* $ so that  the residue vanishes, the constraint $\Sigma-D^*  \succeq 0$ is automatically satisfied.

\section{The Proposed Algorithm}\label{Section_algorithm}
A closed-form solution for Problem \eqref{the_problem0} appears to be out of reach.  However, this Problem appears to be well suited for a coordinate descent type iterative algorithm. 
Such algorithm alternates between solving a minimization problem with respect to $L$ and  a minimization problem with respect to $D$ in the following fashion:
\begin{equation}
\label{alternate_min}
\begin{aligned}
L_k &= \argmin_{L\in\mathcal{L}_{n,r}} \Vert {\Sigma} - L- D_{k-1}\Vert_F\\
D_k &= \argmin_{D\in\mathcal{D}_n} \Vert {\Sigma} - L_{k}- D\Vert_F
\end{aligned}
\end{equation}
 where $L_k$ and $D_k$ denotes the value of $L$ and $D$, respectively, at the $k$-th iteration. Both these sub-problems admit explicit solutions which are provided by the projection operators onto the set $\mathcal{L}_{n,r}$ and $\mathcal{D}_n$, respectively, as described below.
Let $X\in \mathbf{Q}_n$ and consider its spectral decomposition $X=USU^\top$, $U\in\mathbf{O}_n$ and $S=\diag(s_1,\,...,\, s_n)$ with $s_1\geq s_2 \geq ... \geq s_n$ being the eigenvalues of $X$ arranged in decreasing order. Then, the closest   matrix with rank at most $r$ to $X$ in  the Frobenius norm is obtained applying the projector $P_{\mathcal{L}_{n,r}}$:
\begin{equation*}
P_{\mathcal{L}_{n,r}}(X):=U \diag (f_l (s_1), \, ..., \, f_l(s_n) )U^\top
\end{equation*}
with $f_l(\cdot)$ defined as
\begin{equation}
f_l(s_i):=
\begin{cases}
s_i & \text{for } i\leq r \wedge s_i>0 \\
0 & \text{otherwise.} 
\end{cases}
\end{equation}
On the other hand, the projector $P_{\mathcal{D}_n}$ onto the set $\mathcal{D}_n$ is:
\begin{equation}\label{proiezionesud}
P_{\mathcal{D}_n}(X):= \diag(f_d(X_{11}), \, ...., \, f_d(X_{nn}))
\end{equation}
with $f_d(\cdot)$ defined as
\begin{equation}
f_d(X_{ii}):=
\begin{cases}
X_{ii} & \text{if } X_{ii}>0 \\
0 & \text{otherwise.} 
\end{cases}
\end{equation}
Then, at $k$-th iteration the algorithm computes:
\[L_k= P_{\mathcal{L}_{n,r}}({\Sigma}-D_{k-1})\]
\[D_k= P_{\mathcal{D}_n}({\Sigma}-L_k).\]
 The complete procedure is outlined in Algorithm 1: $\varepsilon>0$ is the maximum error allowed in the relative decomposition error, while $N$ represents the maximum number of iterations.
\vspace{0.2cm}
\hrule
\vspace{0.1cm}\noindent
\textbf{Algorithm 1 }
\hrule
\vspace{0.2cm}
\textbf{Input:} ${\Sigma}$, $r$, $\epsilon$, N\\
\indent\textbf{Output:}  $L^*, \;D^*$\\
\indent\textbf{Initialize:} initialize $D$ randomly, i=0\\
\indent\textbf{while} $\Vert {\Sigma} -L-D \Vert_F^2/ \Vert {\Sigma} \Vert_F^2 < \epsilon$ \textbf{and} i $<$ N\\
\indent \indent $L= P_{\mathcal{L}_{n,r}}({\Sigma}-D)$\\
\indent \indent $D= P_{\mathcal{D}_n}({\Sigma}-L)$\\
\indent \indent i=i$+1$\\
\indent\textbf{end while} \\
\indent \indent  $L^*=L$, $D^*=D$\\
\indent\vspace{0.05cm}
\hrule 
\vspace{0.5cm}

\section{Numerical Simulations}\label{sim}

To provide empirical evidence of the convergence properties of the algorithm simulations studies have been performed 
 by using the software Matlab-R2012b on a 2014 laptop MacBook Pro, Quad-i7 2.0 GHz.

 To begin with, we have considered the case of a covariance matrix, $\Sigma$, computed as the sum of a 
randomly generated positive semidefinite low-rank matrix $L$ of dimension $n$ and rank $r$, and a randomly generated positive definite diagonal matrix $D$.
We have performed $200$ Monte Carlo runs with $n=40$ and $r=4$ and  $200$  runs with $n=40$ and $r=10$. The original low-rank and diagonal matrices are recovered with negligible numerical errors. Indeed the following quantities:
\begin{itemize}
\item the  relative decomposition error  on $L+D$: $\Vert \Sigma- L^*- D^* \Vert / \Vert \Sigma\Vert$;
\item the relative error on $L$: $\Vert L-L^*\Vert / \Vert L\Vert$;
\item the relative error on $D$: $\Vert D-D^*\Vert / \Vert D\Vert$;
\end{itemize}
\noindent are all of the order of $10^{-10}$. The average computational time for each experiment is less than five hundredths of a second: in less than half minute all $2\times 200$  runs converged.  


To account for how the algorithm scales with the dimensionality of the problem two further numerical experiments has been conducted.
First, we have considered the case of a fixed rank, $r=8$, and increasing dimensions: $n=20*2^{j}$, with  $j=0\ldots 5$. For each value of $n$, $50$ factor models have been generated and the resulting covariance matrices serve as input for the algorithm. The statistics of the execution time (in seconds) are summarized in Table \ref{table1}. \\
Second, we have considered the case of a fixed $r/n$ ratio of $0.2$, with $n$ taking values $n=20*2^{j}$, with  $j=0\ldots 5$. For each of them, $50$ factor models have been generated and the resulting covariance matrices serve as input for the algorithm. The statistics of the execution time (in seconds) are summarized in Table \ref{table2}. \\
Both experiments provide evidence that the algorithm scales extremely well with dimensionality.

\begin{table}
\centering
\parbox{.45\linewidth}{
\centering
\begin{tabular}{|l| c| c | c|}
\hline
 $n$ & $r$ & mean & st. dev.\\ 
\hline
 20 & 8 & 0.0610 & 0.0545\\
 40 & 8 & 0.0349 & 0.0092\\
 80 & 8 & 0.0642 & 0.0081\\
 160 & 8 & 0.1927 & 0.0173\\
 320 & 8 & 1.2286 & 0.0875\\
 640 & 8 & 5.4973 & 0.2793\\
 1280 & 8 & 26.3725 & 0.9813\\
\hline
\end{tabular}
\caption{\small For each value of $n$ the table displays the mean execution time (in seconds) and standard deviation across $50$ experiments.
}
\label{table1}
}
\quad
\parbox{.45\linewidth}{
\centering
\begin{tabular}{|l| c| c | c|}
\hline
 $n$ & $r$ & mean & st. dev.\\ 
\hline
 20 & 4 & 0.0219 & 0.0452\\
 40 & 8 & 0.0362 & 0.0093\\
 80 & 16 & 0.1069 & 0.0171\\
 160 & 32 & 0.3842 & 0.0444\\
 320 & 64 & 2.7974 & 0.1565\\
 640 & 128 & 14.7181 & 0.6733\\
 1280 & 256 & 85.5031 & 26.8018\\
\hline
\end{tabular}
\caption{\small For each couple of $(n,r)$ the table displays the mean execution time (in seconds) and standard deviation across $50$ experiments.}
\label{table2}
}
\end{table} 

 Finally, we have considered the case of a covariance matrix which admits only approximately a "low-rank plus diagonal" decomposition. This case is of practical interest in factor analysis because typically only an estimate, $\hat{\Sigma}_N$, of $\Sigma$ is available.\\
Given a covariance matrix $\Sigma$ generated as before (which therefore admits an exact "low-rank plus diagonal" decomposition), we have generated a sample of numerosity $N$ from the distribution $\mathcal{N}(0, \Sigma)$ and we have estimated the corresponding sample covariance $\hat{\Sigma}_N$ which serves as input for the algorithm. We have considered the same setting as before with $n=40$, $r=4$ and $n=40$, $r=10$. In both cases for each sample size $N= 200, 500, 1000$ we have performed $200$ Monte Carlo runs. 
The $6\times 200$ simulations took less than 5 minutes to converge  and we observed the following:
\begin{enumerate}
\item
In all the $6\times 200$ simulations the sequence ($D_k,L_k$) produced by Algorithm 1 converged to a stationary point
($L^*,D^*$)
and, as discussed in Proposition \ref{propconv} below, this point is a (at least) local minimum of the cost function.

\item
In all the $6\times 200$ simulations, the inequality $$\Vert L^*+D^* -\hat{\Sigma}_{N} \Vert_F - \Vert\Sigma_{true} -\hat{\Sigma}_{N} \Vert_F\leq 0$$ is satisfied which provides a sanity check on the performance of the proposed algorithm.
In fact, especially for $N=1000$, $\Sigma_{true}$ may be viewed as a good approximation of $\hat{\Sigma}_{N}$ and, on the other hand, we know that, by construction, $\Sigma_{true}$ may be decomposed as the sum of a low rank positive semidefinite matrix and a diagonal positive matrix. Hence, $\Sigma_{true}=L_{true}+D_{true}$ may be viewed as a benchmark which is always outperformed by the decomposition provided by the proposed algorithm.   
\end{enumerate}

\begin{figure}[h!]
	\centering
	\includegraphics[trim={5cm 2cm 4cm 1cm},clip,width=0.5\textwidth]{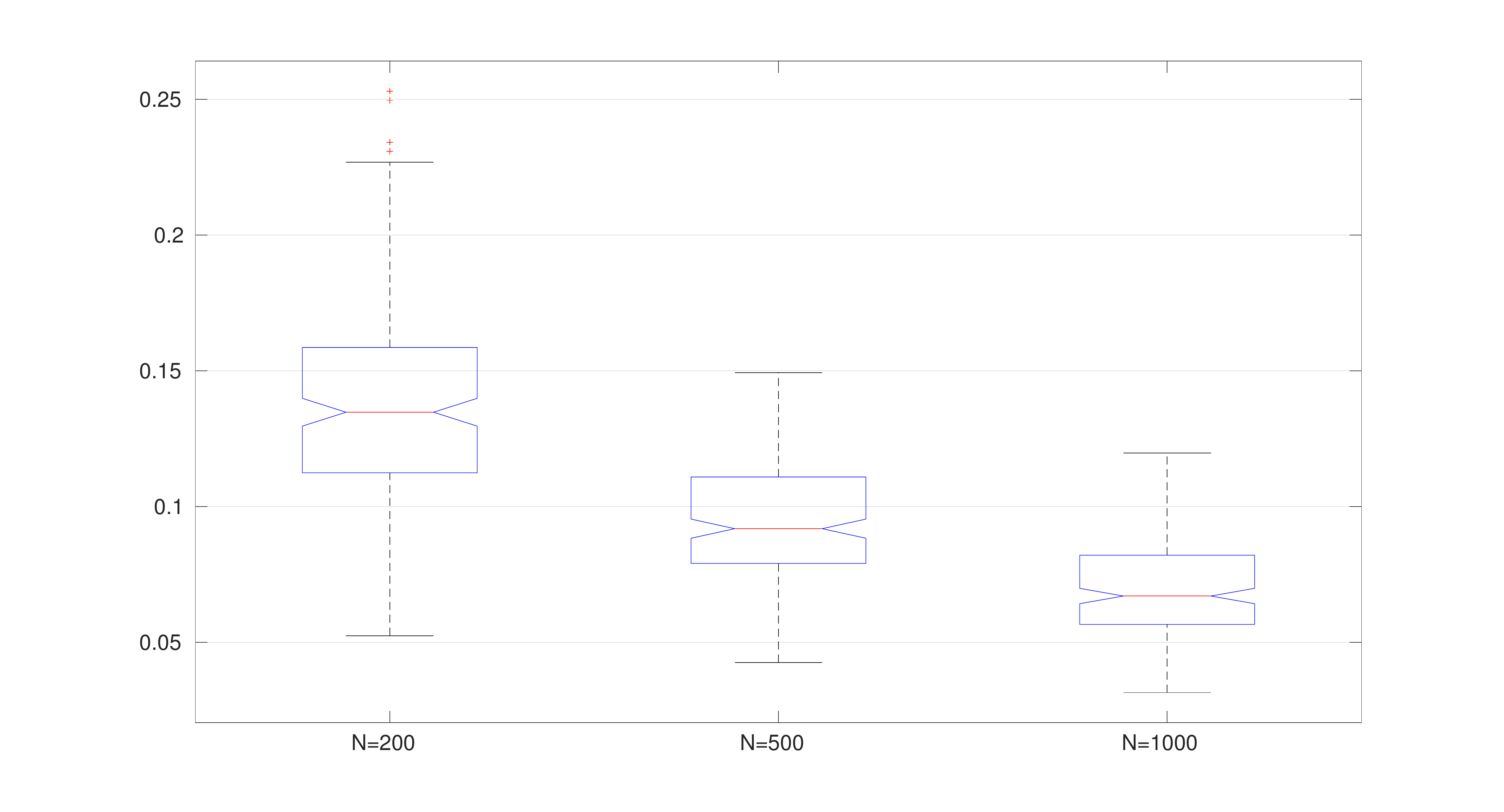}
	\caption{ Case $r=4$. Decomposition errors $\Vert \hat{\Sigma}_N- L^*-D^*\Vert/ \Vert\hat{\Sigma}_N \Vert$ with $N=200$, $N=500$ and $N=1000$. }
\label{fig:figure_5}
\end{figure}

\begin{figure}[h!]
	\centering
	\includegraphics[trim={5cm 2cm 4cm 1cm},clip,width=0.5\textwidth]{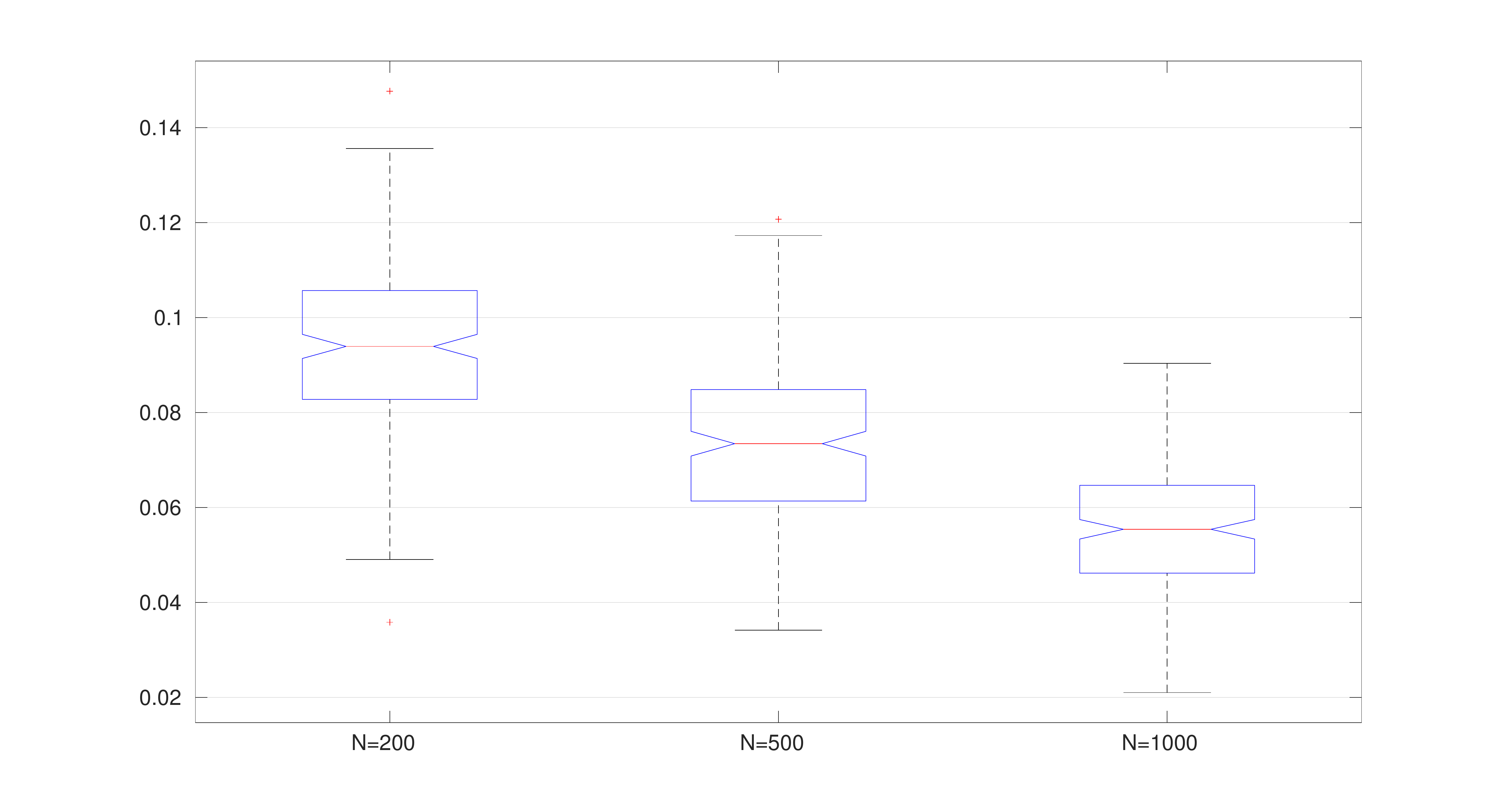}
	\caption{Case $r=10$. Decomposition errors $\Vert \hat{\Sigma}_N- L^*-D^*\Vert/ \Vert \hat{\Sigma}_N \Vert$  with $N=200$, $N=500$ and $N=1000$.}
\label{fig:figure_6}
\end{figure}

\begin{figure}[h!]
	\centering
	\includegraphics[trim={4.5cm 0.5cm 3.5cm 0.5cm},clip,width=0.7\textwidth ,height=11cm]{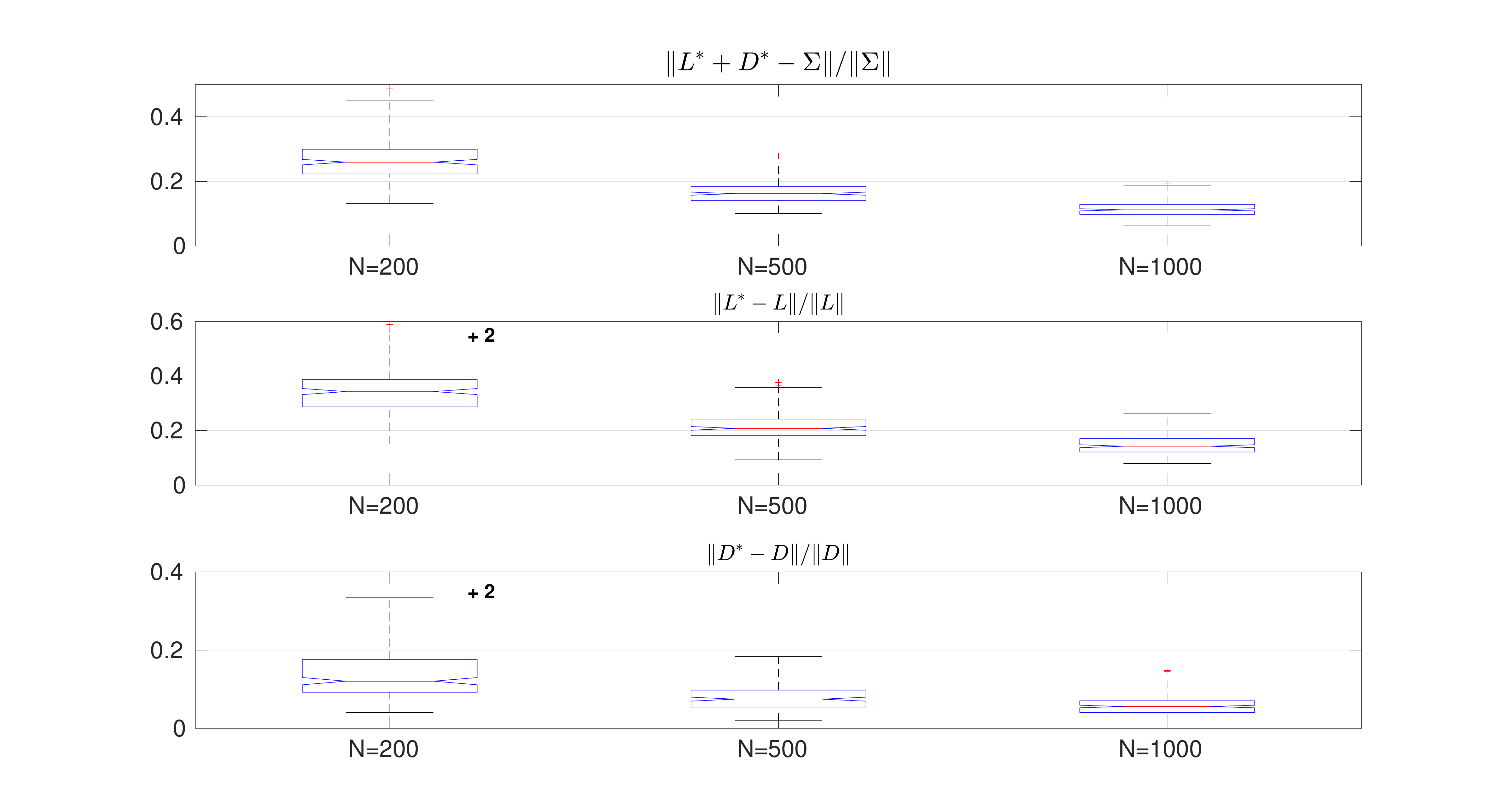}
	\caption{ Case $r=4$. The displayed quantities are: $\Vert \Sigma- L^*-D^*\Vert/ \Vert \Sigma \Vert$, $\Vert L-L^*\Vert/ \Vert L\Vert$ and $\Vert D-D^*\Vert/ \Vert D\Vert$,  where $L^{*}$ and $D^{*}$ represent the estimates with $N=200$, $N=500$ and $N=1000$.}
\label{fig:figure_7}
\end{figure}

\begin{figure}[h!]
	\centering
	\includegraphics[trim={4.5cm 0.5cm 3.5cm 0.5cm},clip,width=0.7\textwidth ,height=11cm]{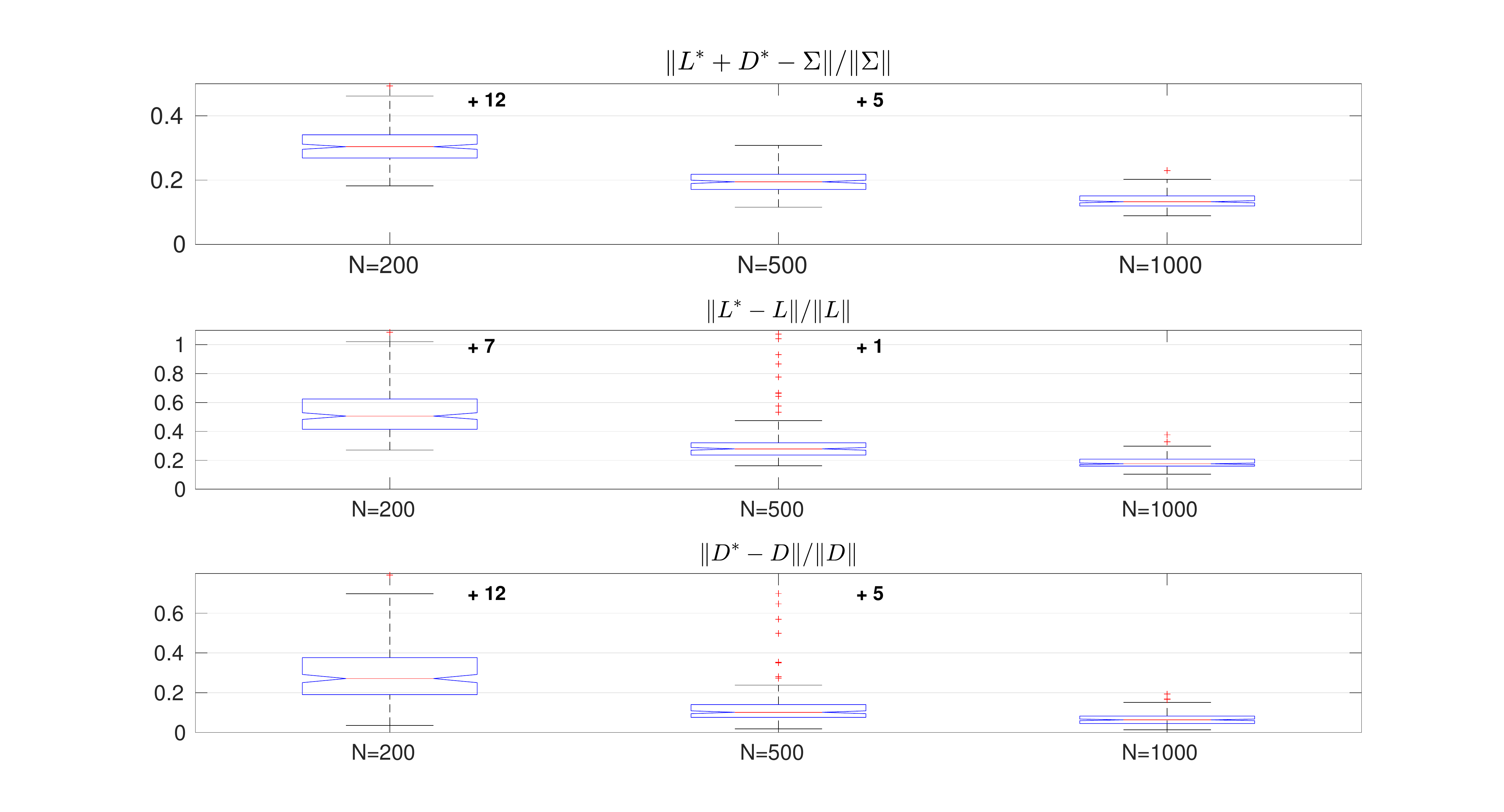}
	\caption{Case $r=10$. The displayed quantities are: $\Vert \Sigma- L^*-D^*\Vert/ \Vert \Sigma \Vert$, $\Vert L-L^*\Vert/ \Vert L\Vert$ and $\Vert D-D^*\Vert/ \Vert D\Vert$, where $L^{*}$ and $D^{*}$ represent the estimates with $N=200$, $N=500$ and $N=1000$. }
\label{fig:figure_8}
\end{figure}

The results for the decomposition error $\Vert \hat{\Sigma}_N- L^*- D^* \Vert / \Vert \hat{\Sigma}\Vert$ are summarized in Figures \ref{fig:figure_5} and \ref{fig:figure_6}. Figures \ref{fig:figure_7} and \ref{fig:figure_8} display the following quantities:
\begin{itemize}
\item the relative  decomposition error on $L+D$: $\Vert \Sigma- L^*-D^*\Vert/ \Vert \Sigma \Vert$; 
\item the relative error on $L$: $\Vert L-L^*\Vert/ \Vert L\Vert$; 
\item the relative error on $D$: $\Vert D-D^*\Vert/ \Vert D\Vert$. \end{itemize}
 
The obtained results appear extremely promising.

\newpage

\vspace{1cm}

\subsection{Application to real data }
\begin{figure}[h!]
	\centering
	\includegraphics[trim={2cm 0cm 3cm 0cm},clip,width=1\textwidth]{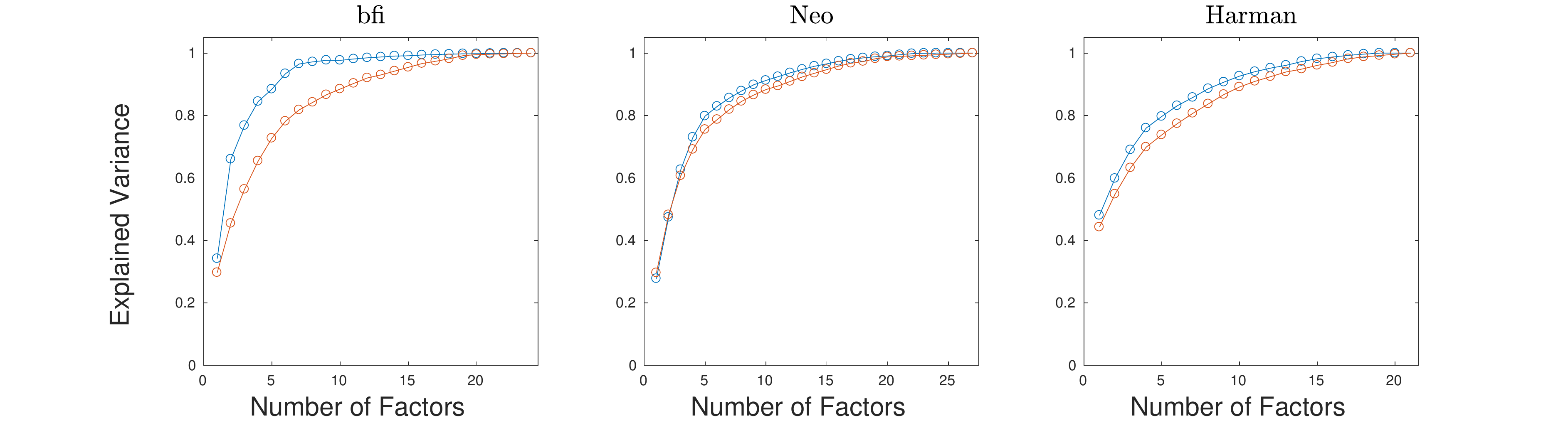}
	\caption{Proportion of variance explained by the proposed method (blue line) and the benchmark method (red line).}
\label{fig:figure_9}
\end{figure}

In this section we investigate the performance of the proposed method on three real world datasets which are popular benchmark in factor analysis:
\begin{itemize}
\item the bfi dataset, from the R library psych, which consists of $2800$ observations on 28 variables: 25 variables represent personality self-reported items and while 3 variables represent demographic variables;
\item the neo dataset, also from the R library psych, which consists of a correlation matrix of size $30 \times 30$ estimated from 1000 observations; 
\item the Harman dataset, from the R library datasets, which consists of a correlation matrix of size $24 \times 24$ estimated from 145 observations: the cross-section represents psychological tests carried out to seventh- and eighth-grade children.
\end{itemize}
These datasets have been used in \cite[Section 5.3]{bertsimas2017certifiably} to compare the performance of their  approach, which minimizes the $q$-norm of the residue (with $q=1$), against different factor analysis methods. This approach can be considered as the state of the art as it outperforms the other available methods.
In this section we take it as benchmark for comparisons and we repeat the analysis in \cite[Section 5.3]{bertsimas2017certifiably}.\\
The adopted measure of performance is the explained variance, defined as $$\sum_{i=1}^r \lambda_i (L^*)/ \sum_{i=1}^n | \lambda_i (\Sigma - D^*)|.$$
For each dataset Problem \eqref{the_problem0} is solved for the values of $r$ considered in \cite{bertsimas2017certifiably}. The results are depicted in Figure \ref{fig:figure_9}.  The proposed method provides a higher amount of explained variance with respect to the method proposed in \cite{bertsimas2017certifiably} that can be considered to be the state-of-the-art as, so far, it outperforms all the available methods. Moreover, our method  shows a flexibility in delivering different models with varying $r$ which is similar to that of the method proposed in \cite{bertsimas2017certifiably}.

In the cases analyzed in these examples the covariance matrix is estimated from a relatively small number of data so that it is reasonable to assume that the residues are not only caused by the presence of less significant latent factors. Thus, we are in the typical situations where our method applies.

\section{Convergence analysis}\label{Section_convergence}
In this section we discuss  the convergence of the proposed algorithm to a local minimum.
First of all we observe that the iterative minimization in \eqref{alternate_min} 
produces a sequence of values  for the objective function that is monotonically non-increasing.
Since the objective  function is clearly bounded from below we have the following obvious result.

\begin{lemm}\label{lemseqfob}
For $h\in\mathbb{N}$, define the sequence  $F_h$ by $F_{h}:=\Vert \Sigma -L_{k}-D_k \Vert_F^2 $ 
for $h=2k$ (even), and $F_{h}:=\Vert \Sigma -L_{k+1}-D_k \Vert_F^2 $, for $h=2k+1$ (odd), where $L_k$, $D_k$ is the sequence produced by Algorithm 1.
Then the sequence $F_h$ is monotonically non-increasing and has limit as $h\rightarrow \infty$.
\end{lemm}

Establishing the convergence for $L_k$ and $D_k$ is less trivial.
We start with $D_k$. To this aim we observe that 
as a consequence of Lemma \ref{lemseqfob}, we have that
 $\varepsilon_k:=F_{2k-1}-F_{2k}$   not only converges to zero but it converges sufficiently fast.

\begin{lemm}\label{lemmaalphaauriga}
Assume that  $\varepsilon_k:=F_{2k-1}-F_{2k}$ tends to zero faster than $1/k^{2q}$ with $q>1$ and
let $D_k$ be the sequence of diagonal matrices produced by Algorithm 1.
Then the sequence $D_k$ converges to a certain diagonal matrix $D^*\in \mathcal{D}_n$.
\end{lemm}
\proof
We have
$$F_{2k}=\Vert \Sigma -L_{k}-D_{k} \Vert_F^2=F_{2k-1}-\varepsilon_k=\Vert \Sigma -L_{k}-D_{k-1} \Vert_F^2-\varepsilon_k.$$ 
Let $s_k(i):=[\Sigma -L_{k}]_{ii}$ be the $i$-th element in the diagonal of $\Sigma -L_{k}$
and $d_k(i):=[D_k]_{ii}$ be the $i$-th element in the diagonal of $D_k$.
Since in (\ref{proiezionesud}) for each $i$, $d_k(i)$ is chosen independently of the others in order to minimize 
$\Vert \Sigma -L_{k}-D_{k} \Vert_F^2$, we have that 
\bea\nn
-\varepsilon_k&=&\Vert \Sigma -L_{k}-D_{k} \Vert_F^2-\Vert \Sigma -L_{k}-D_{k-1} \Vert_F^2\\
\nn
&=&  \sum_{i=1}^n \lbrace [s_k(i)-d_k(i)]^2-[s_k(i)-d_{k-1}(i)]^2 \rbrace \\
\nn
&\leq &
 [s_k(i)-d_k(i)]^2-[s_k(i)-d_{k-1}(i)]^2
 \eea
 which yields
$$\varepsilon_k \geq [d_{k-1}(i)-d_k(i)][d_k(i)+d_{k-1}(i)-2s_k(i)]. $$
Now, we can consider two cases: if $s_k(i)\geq 0$, then the minimizer $d_k(i)$ is equal to 
$s_k(i)$, so that we have
$$\varepsilon_k \geq [d_{k-1}(i)-d_k(i)]^2.$$
If $s_k(i)< 0$, then $d_k(i)=0$  so that we have again
$$\varepsilon_k \geq d_{k-1}(i)[d_{k-1}(i)-2s_k(i)]\geq [d_{k-1}(i)-d_k(i)]^2. $$
In conclusion, in both cases,
we have
$$|d_{k-1}(i)-d_k(i)|\leq \alpha_k:=\sqrt{\varepsilon_k}.$$
As a consequence, we have
\bea\nn
|d_{k+m}(i)-d_k(i)|&\leq & |d_{k+m}(i)-d_{k+m-1}(i)|
+ |d_{k+m-1}(i)-d_{k+m-2}(i)|+\dots \\
\nn
&&
+ |d_{k+1}(i)-d_{k}(i)|\\
\nn
&\leq & \alpha_{k+m}+\dots + \alpha_{k+1}\\
\nn
& = & \sum_{l=1}^m\alpha_{k+l}\\
\nn
&\leq &
\sum_{l=1}^m \frac{M}{(k+l)^{q}}\\
\nn
&\leq &
\sum_{l=1}^\infty \frac{M}{(k+l)^{q}}\\
\nn
&\leq &
 \sum_{h=k+1}^\infty \frac{M}{h^q}
= \left[\sum_{h=1}^\infty \frac{M}{h^q}-\sum_{h=1}^{k+1} \frac{M}{h^q}\right]
\eea
where  $M$ is a constant and $q>1$ so that all the infinite sums converge to a finite value. 
Since we have
\begin{align*}
\lim_{k\rightarrow \infty} \left[\sum_{h=1}^\infty \frac{M}{h^q}-\sum_{h=1}^{k+1} \frac{M}{h^q}\right]=\sum_{h=1}^\infty \frac{M}{h^q}-\lim_{k\rightarrow \infty} \sum_{h=1}^{k+1} \frac{M}{h^q}=0,
\end{align*} 
we can conclude that $\lim_{l,k\rightarrow\infty}|d_l(i)-d_k(i)|=0$, so that $d_k(i)$ is a Cauchy sequence and hence  it converges.
Since this holds for each $i=1,\dots,n$, we have that the sequence $D_k$ converges to a certain diagonal matrix $D^*$.
Finally, since $\mathcal{D}_n$ is closed, clearly $D^*\in \mathcal{D}_n$.
\qed

For the convergence of the sequence $L_k$ we need to rule out a pathological situation.

\begin{lemm}\label{lconv}
Under the assumptions of Lemma \ref{lemmaalphaauriga}, let $S:=\Sigma-D^*$ with $D^*$ being the limit of the sequence of diagonal matrices $D_k$ produced by Algorithm 1.
If $S$ has $n$ distinct eigenvalues then
the sequence of rank $r$ matrices $L_k$ produced by Algorithm 1 converges to a rank $r$
matrix $L^*$. 
\end{lemm}
\proof
Let $s_1> s_2 > ... > s_n$ be the eigenvalues of $S$ arranged in decreasing order.
By continuity of the eigenvalues, for a sufficiently large $k$, $S_k:=\Sigma-D_k$
has $n$ distinct eigenvalues $s_{k,1}> s_{k,2} > \dots > s_{k,n}$ and 
$\lim_{k\rightarrow\infty}s_{k,i}=s_i$.
 
According to \cite[Chapter 9, Theorem 8]{lax07}, for each $i=1,\dots,n$, we can select an eigenvector
(and hence a normalized eigenvector) $v_{k,i}$ of $S_k$ associated with the eigenvalue 
$s_{k,i}$ in such a way that $v_{k,i}$ converges to a normalized eigenvector of $S$ associated with the eigenvalue $s_i$.
Now recall that
$$
L_{k+1}=P_{\mathcal{L}_{n,r}}(S_k)=U_k \diag (f_l (s_{k,1}), \, ..., \, f_l(s_{k,n}) )U_k^\top
$$
where the $i$-th column of $U_k$ is a normalized eigenvector of $S_k$ associated with the eigenvalue $s_{k,i}$.
As a normalized eigenvector is unique up to its sign, we have
$U_k=V_k\Delta_k$ with  $V_k:=[v_{k,1}\mid v_{k,2} \mid\dots\mid v_{k,n}]$
and $\Delta_k$ is a diagonal matrix whose diagonal entries can only be $\pm 1$.
We easily see that  the contribution of the $\Delta_k$  cancels and we have
$$
L_{k+1}=V_k \diag (f_l (s_{k,1}), \, ..., \, f_l(s_{k,n}) )V_k^\top
$$
so that $L_{k+1}$ is given by the product of three matrices each one of which converges as $k$ tends to infinity.
\qed

\begin{propo}\label{propconv}
Assume that the hypothesis of Lemma \ref{lconv} holds and that the matrix $L^*$ defined in the same lemma has rank $r$. Assume also that the tangent space of $\mathcal{L}_{n,r}$ at 
$L^*$ does not contain diagonal matrices. Then the sequence ($D_k,L_k$) produced by Algorithm 1 converges to a point corresponding to a  local minimum of the cost function.
\end{propo}
\proof
By the previous results, we know that $D_k$ converges to $D^*$ and $L_k$ converges to $L^*$.
Assume by contradiction that ($D^*,L^*$) is not a  minimum.
Then,  for  any $\varepsilon>0$, there exists $\delta D$ and $\delta L$ such that
$\|\delta D\|_F<\varepsilon$, $\|\delta L\|_F<\varepsilon$, $(L+\delta L)\in {\mathcal{L}_{n,r}}$,
$(D+\delta D)\in {\mathcal{D}_n}$ and
$$
\|\Sigma -L^*-D^* \|_F^2> \|\Sigma -L^*-D^* -\delta L- \delta D \|_F^2.
$$
Now let $\delta T$ be the projection of $\delta L$ on the tangent space of $\mathcal{L}_{n,r}$ at 
$L^*$.  For a sufficiently small  $\varepsilon$ we have
$$
\|\Sigma -L^*-D^* \|_F^2\geq \|\Sigma -L^*-D^* -\delta T- \delta D \|_F^2.
$$

By setting $R:=\Sigma -L^*-D^*$ and computing the Frobenius norms in the previous formula, we get
$$
2(\tr[R\delta T] + \tr[R\delta D]) - \|\delta T+ \delta D \|_F^2\geq 0.
$$
By assumption $\delta T+ \delta D\neq 0$ so that 
at least one of the two quantities  $2\tr[R\delta L]$ and  $2\tr[R\delta D]$ is positive.
In the case of $\tr[R\delta D]>0$ we have that for all $\kappa$ sufficiently small,
\bea\nn
\min_{D\in\mathcal{D}_n} \Vert {\Sigma} - L^*- D\Vert_F^2 &\leq& \Vert {\Sigma} - L^*- D^*-\kappa \delta D\Vert_F^2\\
\nn
&=& \Vert R\Vert_F^2 + \kappa^2 \Vert \delta D\Vert_F^2- 2\kappa \tr[R\delta D] \\
\nn
&<&\Vert R\Vert_F^2
\eea
which is contradiction because we know that the algorithm converged so that 
$\min_{D\in\mathcal{D}_n} \Vert {\Sigma}- L^*- D\Vert_F^2=\Vert R\Vert_F^2$.

In the case of $\tr[R\delta T]>0$ we have that
\beq
\min_{L\in\mathcal{L}_{n,r}} \Vert {\Sigma} - L- D^*\Vert_F^2 \leq  \Vert {\Sigma} - D^* -P_{\mathcal{L}_{n,r}}(L^*+\kappa\delta T )  \Vert_F^2
\eeq
where $P_{\mathcal{L}_{n,r}}(\cdot)$ is the projection onto $\mathcal{L}_{n,r}$.
Thus we have
$P_{\mathcal{L}_{n,r}}(L^*+\kappa\delta T ) =L^*+\kappa\delta T + E$
where $\lim_{\kappa\rightarrow 0} \Vert E\Vert_F/\kappa =0$.
 
Thus, for $\kappa>0$ sufficiently small, we have
\bea\nn 
q&:=& \Vert {\Sigma} - D^* -P_{\mathcal{L}_{n,r}}(L^*+\kappa\delta T )  \Vert_F^2 \\
\nn
&=&\Vert {\Sigma} - D^* -L^* - \kappa\delta T - E  \Vert_F^2
\\
\nn
&=&\Vert R  \Vert_F^2+ \kappa^2 \Vert \delta T  \Vert_F^2 +   \Vert E  \Vert_F^2
 -2 \kappa\tr(R\delta T) -2\tr(R E) +2\kappa\tr(\delta T E)\\
\nn 
&<&\Vert R  \Vert_F^2.
\eea
In conclusion, we have
\beq
\min_{L\in\mathcal{L}_{n,r}} \Vert {\Sigma} - L- D^*\Vert_F^2 < \Vert R  \Vert_F^2,
\eeq
that, as in the previous case leads to a contradiction.
\qed

\begin{rem}
We believe that the assumption of Lemma \ref{lconv} can be weakened  that the results still hold
if $s_r > s_{r+1}$ where $s_1\geq s_2 \geq ... \geq s_n$ are the eigenvalues of $S$ repeated according to their algebraic multiplicity and arranged in decreasing order.
The proof of this fact seems, however, very delicate because of some issues on the continuity of eigenspaces under small perturbations.
\end{rem}

\begin{rem}
It is quite intuitive that the conditions of Proposition \ref{propconv} are not very stringent: in fact in all the practical situations that we have studied in simulations those  conditions are satisfied.
\end{rem}

\section{An Alternating Projection Type Algorithm}\label{pojection-approach}

In this section we present our algorithm under a different perspective that may be useful in addressing questions on the properties of the proposed method.
In fact, by suitably translating ${\mathcal{D}_n}$, we  easily see that this method can be viewed as an alternating projection type algorithm for which a very rich literature has been developed.
To this aim,
define
\beq    
\begin{aligned}
\mathcal{\tilde{D}}_n:=\Sigma-\mathcal{D}_n=\lbrace X\in\mathbf{Q}_n: \ofd(X)=\ofd({\Sigma}), \, X_{ii}\leq {\Sigma}_{ii}, i=1,...,n \rbrace
\end{aligned}
\eeq
and notice that the projection in this affine set is easily obtained as:
\begin{equation}
P_{\mathcal{\tilde{D}}_n}(X):=
\begin{cases}
X_{ij}=X_{ij} & \text{for } i= j \wedge X_{ii}<{\Sigma}_{ii} \\
X_{ij}={\Sigma}_{ij} & \text{for } (i= j \wedge X_{ii}\geq {\Sigma}_{ii}) \vee i\neq j. 
\end{cases}
\end{equation}
We consider now the sequences $L_k$ and $D_k$ produced by our algorithm.
We recall that our  $D_k$  is given by $D_k = P_{\mathcal{D}_n}({\Sigma}-L_k)$.
By taking this formula into account, a direct computation shows that the matrix
$\tilde{D}_k:=P_{\mathcal{\tilde{D}}_n}(L_k)$ may be written as
$\Sigma-D_k$ so that, in view of the formula $L_k= P_{\mathcal{L}_{n,r}}({\Sigma}-D_{k-1})$,
we immediately get that
\[L_{k+1}=P_{\mathcal{L}_{n,r}}(P_{\mathcal{\tilde{D}}_n}(L_k))\]
which shows that the iteration for $L_k$ is the result of  an alternating projection algorithm.
These kind of algorithms burst a long tradition which dates back to Von Neumann in the '30s. 
While for alternating projection onto convex sets the convergence results are well established, for the non-convex case much less is known.
In our case $\mathcal{\tilde{D}}_n$ is a convex set of dimension $n$, but the set $\mathcal{L}_{n,r}$ is a non-convex embedded manifold of $\mathbb{R}^{n\times n}$ with dimension $nr-r(r-1)/2$ and it is smooth  at those points for which the rank is exactly $r$. 
In \cite{lewis2008alternating} a proof of local convergence (at a linear rate) for alternating projection onto smooth manifolds is provided  under the assumption of transversal intersection.
In our case, transversal intersection cannot hold when $r$ is small with respect to $n$ but it may be possible to generalise that approach to provide a further analysis of the algorithm properties and, in particular, of its convergence rate.

Finally, the set $\tilde{D}_n$ is particularly interesting because of the following interpretation that is particularly  evident when $r$ is such that $\Sigma$ can be decomposed exactly as  $L^*+D^*$ so that $(L^*,D^*)$ is clearly an optimal solution of \eqref{the_problem0}. In this case, $D^*=\Sigma-L^*$ and thus $\Sigma-L^*\in\mathcal D_n$. The latter condition is equivalent to the condition $L^*\in\tilde{\mathcal D_n}$
Therefore the problem \eqref{the_problem0} can reformulated only in terms of $L$ as follows:
\begin{equation}
\begin{aligned}
\label{the_problem}
L^*:=\argmin_{L\in \mathcal{L}_{n,r} \cap \mathcal{\tilde{D}}_n,} \quad & \Vert {\Sigma} -L\Vert_F^2 
\end{aligned}
\end{equation}
It is worth noting that the objective function in \eqref{the_problem} does not take into account the covariance matrix of the idiosyncratic noise, i.e. such a matrix is understood as the covariance matrix of a noise random vector. The latter is in the same spirit of  \cite{stock1998diffusion} wherein the factor loading matrix $A$ is given by solving a least squares problem for the linear regression model \eqref{fact_model} and the idiosyncratic component is treated as noise.

\section{Conclusions}
 We have proposed an alternating minimization algorithm for decomposing a covariance matrix as sum of a low rank matrix, whose maximal rank is a priori fixed, plus a diagonal matrix. The latter minimizes the residue among the covariance matrix and the additive decomposition. Simulation results showed that the algorithm performs extremely well and converges very rapidly to the solution. Finally, we have proved that, under reasonable assumptions, such algorithm converges to a solution which is a local minimum for the residue.

\vspace{\fill}\pagebreak



\begin{thebibliography}{10}

\bibitem{agarwal2012noisy}
A.~Agarwal, S.~Negahban, and M.~J. Wainwright.
\newblock Noisy matrix decomposition via convex relaxation: Optimal rates in
  high dimensions.
\newblock {\em The Annals of Statistics}, pages 1171--1197, 2012.

\bibitem{bai2012statistical}
J.~Bai, K.~Li, et~al.
\newblock Statistical analysis of factor models of high dimension.
\newblock {\em The Annals of Statistics}, 40(1):436--465, 2012.

\bibitem{bai2002determining}
J.~Bai and S.~Ng.
\newblock Determining the number of factors in approximate factor models.
\newblock {\em Econometrica}, 70(1):191--221, 2002.

\bibitem{bai2008large}
J.~Bai, S.~Ng, et~al.
\newblock Large dimensional factor analysis.
\newblock {\em Foundations and Trends{\textregistered} in Econometrics},
  3(2):89--163, 2008.

\bibitem{Bekker-deLeeuw_1987}
P.~A. Bekker and J.~de~Leeuw.
\newblock The rank of reduced dispersion matrices.
\newblock {\em Psychometrika}, 52(1):125?--135, 1987.

\bibitem{bertsimas2017certifiably}
D.~Bertsimas, M.~S. Copenhaver, and R.~Mazumder.
\newblock Certifiably optimal low rank factor analysis.
\newblock {\em Journal of Machine Learning Research}, 18(29):1--53, 2017.

\bibitem{bottegal2015modeling}
G.~Bottegal and G.~Picci.
\newblock Modeling complex systems by generalized factor analysis.
\newblock {\em IEEE Transactions on Automatic Control}, 60(3):759--774, 2015.

\bibitem{BURT_1909}
C.~Burt.
\newblock Experimental tests of general intelligence.
\newblock {\em British Journal of Psychology, 1904-1920}, 3(1/2):94--177, 1909.

\bibitem{ciccone2017factor}
V.~Ciccone, A.~Ferrante, and M.~Zorzi.
\newblock Is factor analysis viable for real world data?
\newblock {\em arXiv preprint arXiv:1709.01168}, 2017.

\bibitem{deistler2007}
M.~Deistler and C.~Zinner.
\newblock Modelling high-dimensional time series by generalized linear dynamic
  factor models: An introductory survey.
\newblock {\em Communications in Information \& Systems}, 7(2):153--166, 2007.

\bibitem{MIN_RANK_SHAPIRO_1982}
G.~Della~Riccia and A.~Shapiro.
\newblock Minimum rank and minimum trace of covariance matrices.
\newblock {\em Psychometrika}, 47:443--448, 1982.

\bibitem{fan2013large}
J.~Fan, Y.~Liao, and M.~Mincheva.
\newblock Large covariance estimation by thresholding principal orthogonal
  complements.
\newblock {\em Journal of the Royal Statistical Society: Series B (Statistical
  Methodology)}, 75(4):603--680, 2013.

\bibitem{FAZEL_MIN_RANK_APPLICATIONS_2002}
M.~Fazel.
\newblock Matrix rank minimization with applications.
\newblock {\em Elec. Eng. Dept. Stanford University}, 54:1--130, 2002.

\bibitem{FAZEL_MINIMUM_RANK_2004}
M.~Fazel, H.~Hindi, and S.~Boyd.
\newblock Rank minimization and applications in system theory.
\newblock In {\em Proceedings of the American Control Conference}, volume~4,
  pages 3273--3278, Jun. 2004.

\bibitem{finesso2016factor}
L.~Finesso and P.~Spreij.
\newblock Factor analysis models via i-divergence optimization.
\newblock {\em psychometrika}, 81(3):702--726, 2016.

\bibitem{GEWEKE_DYNAMIC}
J.~Geweke.
\newblock The dynamic factor analysis of economic time series models.
\newblock In {\em Latent Variables in Socio-Economic Models}, SSRI workshop
  series, pages 365--383. North-Holland, 1977.

\bibitem{guttman1954some}
L.~Guttman.
\newblock Some necessary conditions for common-factor analysis.
\newblock {\em Psychometrika}, 19(2):149--161, 1954.

\bibitem{harman1966factor}
H.~H. Harman and W.~H. Jones.
\newblock Factor analysis by minimizing residuals (minres).
\newblock {\em Psychometrika}, 31(3):351--368, 1966.

\bibitem{DEISTLER_1997}
C.~Heij, W.~Scherrer, and M.~Deistler.
\newblock System identification by dynamic factor models.
\newblock {\em SIAM Journal on Control and Optimization}, 35(6):1924--1951,
  1997.

\bibitem{kaiser1958varimax}
Henry~F Kaiser.
\newblock The varimax criterion for analytic rotation in factor analysis.
\newblock {\em Psychometrika}, 23(3):187--200, 1958.

\bibitem{lam2012factor}
C.~Lam, Q.~Yao, et~al.
\newblock Factor modeling for high-dimensional time series: inference for the
  number of factors.
\newblock {\em The Annals of Statistics}, 40(2):694--726, 2012.

\bibitem{lax07}
P.~D. Lax.
\newblock {\em Linear Algebra and Its Applications}.
\newblock Wiley-Interscience, second edition, 2007.

\bibitem{lewis2008alternating}
A.~S. Lewis and J.~Malick.
\newblock Alternating projections on manifolds.
\newblock {\em Mathematics of Operations Research}, 33(1):216--234, 2008.

\bibitem{ning2015linear}
L.~Ning, T.~T. Georgiou, A.~Tannenbaum, and S.~P. Boyd.
\newblock Linear models based on noisy data and the {F}risch scheme.
\newblock {\em SIAM Review}, 57(2):167--197, 2015.

\bibitem{picci1986dynamic}
G.~Picci and S.~Pinzoni.
\newblock Dynamic factor-analysis models for stationary processes.
\newblock {\em IMA Journal of Mathematical Control and Information},
  3(2-3):185--210, 1986.

\bibitem{reiersol1950identifiability}
O.~Reiers{\o}l.
\newblock Identifiability of a linear relation between variables which are
  subject to error.
\newblock {\em Econometrica: Journal of the Econometric Society}, pages
  375--389, 1950.

\bibitem{scherrer1998structure}
W.~Scherrer and M.~Deistler.
\newblock A structure theory for linear dynamic errors-in-variables models.
\newblock {\em SIAM Journal on Control and Optimization}, 36(6):2148--2175,
  1998.

\bibitem{shapiro1982rank}
A.~Shapiro.
\newblock Rank-reducibility of a symmetric matrix and sampling theory of
  minimum trace factor analysis.
\newblock {\em Psychometrika}, 47(2):187--199, 1982.

\bibitem{shapiro2002statistical}
A.~Shapiro and J.~MF Ten~Berge.
\newblock Statistical inference of minimum rank factor analysis.
\newblock {\em Psychometrika}, 67(1):79--94, 2002.

\bibitem{spearman_1904}
C.~Spearman.
\newblock {"General Intelligence," Objectively Determined and Measured}.
\newblock {\em American Journal of Psychology}, 15:201--293, 1904.

\bibitem{Spearman-Holzinger-24}
C.~Spearman and K.~J. Holzinger.
\newblock The sampling error in the theory of two factor.
\newblock {\em British Journal of Psychology}, 15:17--19, 1924.

\bibitem{stock1998diffusion}
J.~H. Stock and M.~W. Watson.
\newblock Diffusion indexes.
\newblock Technical report, National bureau of economic research, 1998.

\bibitem{tucker1973reliability}
L.~R. Tucker and C.~Lewis.
\newblock A reliability coefficient for maximum likelihood factor analysis.
\newblock {\em Psychometrika}, 38(1):1--10, 1973.

\bibitem{zorzi2016ar}
M.~Zorzi and R.~Sepulchre.
\newblock Ar identification of latent-variable graphical models.
\newblock {\em IEEE Transactions on Automatic Control}, 61(9):2327--2340, 2016.

\end{thebibliography}
\end{document}